\documentclass[12pt,reqno]{amsart}
\usepackage{amssymb,latexsym}
\newtheorem{theorem}{Theorem}
\newtheorem{lemma}{Lemma}
\newtheorem*{cortheorem}{Corona Theorem}
\newtheorem*{theoremAM}{Theorem (Agler and McCarthy)}
\newtheorem*{theoremAmar}{Theorem (Amar)}

\newcommand{\wlim}{\operatorname{weak\,limit}}

\begin{document}
\title[Toeplitz Corona Theorems]
{Toeplitz Corona Theorems for the Polydisk and the Unit Ball}
\author[T. T. Trent]{Tavan T. Trent}
\address{Department of Mathematics\\
         The University of Alabama\\
         Box 870350\\
         Tuscaloosa, AL  USA 35487-0350}
\email{ttrent@as.ua.edu}
\thanks{The research of the first author was partially supported by the Fields Institute.}
\author[B. D. Wick]{Brett D. Wick}
\address{Department of Mathematics\\ University of South Carolina\\ LeConte College\\ 1523 Greene Street\\ Columbia, SC USA 29206}
\address{Present:  The Fields Institute\\ 222 College Street, 2nd Floor\\ Toronto, Ontario\\ M5T 3J1 Canada}
\email{wick@math.sc.edu}
\thanks{The research of the second author was partially supported by the Fields Institute and by NSF DMS Grant \# 0752703.}
\begin{abstract}
The main purpose of this paper is to extend and refine some work of Agler-McCarthy and Amar concerning the Corona problem for the polydisk and the unit ball in $\mathbb{C}^n$.
\end{abstract}
\maketitle

\medskip
The main purpose of this paper is to extend and refine some work of Agler-McCarthy \cite{AM} and Amar \cite{A} concerning the Corona problem for the polydisk and the unit ball in $\mathbb{C}^n$.  In 1962, Carleson \cite{C} proved his famous Corona theorem with bounds:

\begin{cortheorem}
Let $\{ f_j\}_{j=1}^m \subseteq H^{\infty}(D)$.  Assume that
\[
0 < \epsilon^2 \le \sum_{j=1}^m \, |f_j(z)|^2 \le 1 \; \, \mbox{for all } \, z \in D.
\]
There exists a constant  $C(\epsilon,m)< \infty$  and  $\{ g_j\}_{j=1}^m \subseteq H^{\infty}(D)$, so that
\[
\sup_{z \in D} \, \sum_{j=1}^m \, |g_j(z)|^2 \le C(\epsilon,m)^2 \; \mbox{and} \; \sum_{j=1}^m \, f_j(z)g_j(z) = 1 \; \, \mbox{for all } \, z \in D.
\]
\end{cortheorem}

\noindent
The Corona theorem and especially the techniques utilized in its proof have been very influential.  See, for example, Garnett \cite{Ga}.  Among many questions raised by this theorem, we wish to consider the analogous Corona problem for the polydisk and the unit ball in $\mathbb{C}^n$.

\medskip
We will need some notation:
\begin{align}
& D^n &\;  & \{z = (z_1, \dots, z_n) \in \mathbb{C}^n : \, |z_j| < 1 \;\; \mbox{for } \, j=1, \dots, n\} \notag \\
& B^n & \; &\{z = (z_1, \dots, z_n) \in \mathbb{C}^n : \; \sum_{j=1}^n \, |z_j|^2 < 1\} \notag \\
& \Omega &\; &\mbox{denotes either } \, D^n \; \mbox{or } \, B^n \notag \\
& \partial \Omega &\; &\mbox{denotes } \; \partial B^n \; \mbox{if } \; \Omega = B^n \; \mbox{or else} \, , \mbox{if } \, \Omega = D^n, \notag \\
& &{ } &\partial \Omega \;\; \mbox{denotes } \; T^n, \mbox{the \underline{distinguished} boundary of } \, D^n \notag \\
& \sigma &\; &\mbox{denotes normalized Lebesgue measure on } \, \partial \Omega \notag \\
& \mu &\;  &\mbox{denotes a probability measure on } \, \partial \Omega \notag \\
& P^2(\mu) &\; &\{ \mbox{analytic polynomials on } \mathbb{C}^n \}^{-L^2(\mu)} \notag \\
& F(z) &\; &(f_1(z), f_2(z), \dots) \;\; \mbox{for } \,\; f_j \in H^{\infty}(\Omega) \notag \\
& \mathcal{F}(z) &\; &[f_{ij}(z)]_{i,j = 1}^{\infty} \; \mbox{for} \; f_{ij} \in H^{\infty}(\Omega) \notag \\
& T_F^{\mu} &\; &\mbox{the multiplication operator from } \, \underset{1}{\overset{\infty}{\oplus}} \, P^2(\mu) \; \mbox{into } \, P^2(\mu) \; \mbox{by} \, F\notag \\
& T_{\mathcal{F}}^{\mu} &\; &\mbox{the multiplication operator from }  \underset{1}{\overset{\infty}{\oplus}} \, P^2(\mu) \; \mbox{into }  \underset{1}{\overset{\infty}{\oplus}} \, P^2(\mu) \; \mbox{by} \, \mathcal{F}\notag \\
& I_{\mu} &\; &\mbox{the identity operator on } \, P^2(\mu) \notag \\
& \mathcal{H} &\; &\{H \in H^{\infty}(\Omega): \, H \; \mbox{nonvanishing in }  \Omega, \, \frac 1{H} \in L^{\infty}(\partial \Omega, d\sigma), \notag \\
& &\; &\mbox{ and }  \Vert H\Vert_2 = 1 \} \notag \\
& T_F^H &\; &\mbox{defined as above, when} \; d\mu = |H|^2 d\sigma\; \mbox{and }\, H \in \mathcal{H} \notag
\end{align}

\bigskip
For the case of the bidisk, $D^2$, Agler and McCarthy proved the following:
\begin{theoremAM}
Let  $\{ f_j\}_{j=1}^m \subseteq H^{\infty}(D^2)$.  Then there exist  $ \{ g_j\}_{j=1}^m \subseteq H^{\infty}(D^2)$  with
\[
\sum_{j=1}^m \, f_jg_j \equiv 1 \;\; \mbox{and } \; \sup_{z \in D^2} \, \sum_{j=1}^m \, |g_j(z)|^2 \le \frac 1{\delta^2}
\]
if and only if
\[
T_F^{\mu}(T_F^{\mu})^* \ge \delta^2 I_{\mu}
\]
for all probability measures $\mu$ on $T^2$.
\end{theoremAM}

\medskip
Although the Agler-McCarthy theorem and its proof seemed to be restricted to $n=2$ by the classical and beautiful counterexample of Parrot \cite{P}; nevertheless, Amar managed to extend it to $D^n$ (and to $B^n$).
\begin{theoremAmar}
Let  $\{ f_j\}_{j=1}^m \subseteq H^{\infty}(\Omega)$.  Then there exist  $ \{ g_j\}_{j=1}^m \subseteq H^{\infty}(\Omega)$  with
\[
\sum_{j=1}^m \, f_jg_j \equiv 1 \;\; \mbox{and } \; \sup_{z \in \Omega} \, \sum_{j=1}^m \, |g_j(z)|^2 \le \frac 1{\delta^2}
\]
if and only if
\[
T_F^{\mu}(T_F^{\mu})^* \ge \delta^2 I_{\mu}
\]
for all probability measures $\mu$ on $\partial \Omega$.
\end{theoremAmar}

\medskip
In other words, Amar shows that for $\{ f_j\}_{j=1}^m \subseteq H^{\infty}(\Omega)$ and $\delta >0$ the following are equivalent:

\medskip
(i) There exist  $ \{ g_j\}_{j=1}^m \subseteq H^{\infty}(\Omega)$  with
\[
\sum_{j=1}^m \, f_jg_j = 1 \;\;\text{on} \,\, \Omega \,\, \mbox{and } \; \sup_{z \in \Omega} \, \sum_{j=1}^m \, |g_j(z)|^2 \le \frac 1{\delta^2}. 
\]

(ii) For all probability measures $\mu$ on $ \partial \Omega$ and all $h \in P^2(\mu)$ there exists $ \{ k_j\}_{j=1}^m \subseteq H^{\infty}(\Omega)$
\[
\sum_{j=1}^m \, f_jk_j = h \;\; \mbox{and } \;  \, \sum_{j=1}^m \, ||k_j||_{\mu}^2 \le \frac 1{\delta^2}||h||_{\mu}^2 
\]

\bigskip
By results of Andersson-Carlsson \cite{AC} for the unit ball and Varopoulos \cite{V}, Li \cite{L}, Lin \cite{Li}, Trent \cite{Tr}, and Treil-Wick \cite{TW} for the polydisk case, we know that if the input functions are bounded away from $0$ on $\Omega$, we have an  $H^p(\Omega)$  Corona theorem for  $1 \le p < \infty$.  That is, if
\[
\{ f_j\}_{j=1}^{\infty} \subseteq H^{\infty}(\Omega) \;\; \mbox{and } \; 0 < \epsilon^2 \le \underset{z \in \Omega}{\mathrm{inf}} \, \sum_{j=1}^{\infty} \, |f_j(z)|^2 \le 1,
\]
then for  $1 \le p < \infty$  there exists a  $\delta_p > 0$  so that
\[
T_FT_F^* \ge \delta_p^2 I_{H^p(\Omega)},
\]
where \, $F = (f_1, f_2, \dots)$. Unfortunately, the best of these estimates have $ \delta_p \downarrow 0\; \mbox{as} \; p \uparrow \infty$.

\bigskip
Thus Amar's theorem tells us that a solution to the Corona problem for  $H^{\infty}(\Omega)$  follows from the following statement:
$$
T_FT_F^* \ge \delta^2 I, \; \mbox{for} \; \delta > 0 \; \overset{\mbox{?}}{\Rightarrow} \;  \exists \; \epsilon > 0 \;\, \mbox{ such that } T_p^{\mu}(T_p^{\mu})^* \ge \epsilon^2 I_{\mu} \notag
$$
for all probability measures $\mu$ on $ \partial \Omega$.

\medskip
Of course, necessity in Amar's theorem is trivial; so we will concentrate on weakening the  sufficient conditions to get the same Corona output.

\medskip
We will extend Amar's theorem to an infinite number of input functions and refine his theorem, so that we need only consider probability measures, $\mu$, of the form  $|H|^2 \, d \sigma$, where  $H \in \mathcal{H}$.  In addition, we weaken the hypotheses to just have our operators dominate a certain rank one operator.  We begin with a series of lemmas.

\medskip
\begin{lemma}
Let  $\mathcal{F}(z) = [f_{ij}(z)]_{i,j = 1}^{\infty}$, $f_{ij} \in H^{\infty}(\Omega)$.
Then
\[
\Vert T_{\mathcal{F}} \Vert_{B(\underset{1}{\overset{\infty}{\oplus}}\, H^2(\Omega))} = \, \sup_{z \in \Omega} \, \Vert \mathcal{F}(z)\Vert_{B(l^2)}.
\]
\end{lemma}

\medskip

\begin{proof}
Let $ \underline h \in \underset{1}{\overset{\infty}{\oplus}}\,H^2(\Omega)$.  Then
\begin{align}
\Vert T_{\mathcal{F}} \underline h\Vert^2_{\underset{1}{\overset{\infty}{\oplus}}\, H^2(\Omega)} & = \sup_{0\leq r < 1} \left(\int_{\partial \Omega} \, \Vert \mathcal{F}(re^{it}) \underline h(re^{it}) \Vert _{l^2}^2 \, d \sigma\right) \notag \\
& \le \sup_{z \in \Omega} \Vert \mathcal{F}(z) \Vert_{B(l^2)}^2 \, \sup_{0\leq r < 1} \left(\int_{\partial \Omega} \,  \, \Vert \underline h (re^{it}) \Vert_{l^2}^2 \, d\sigma\right) \notag \\
& \le \sup_{z \in \Omega} \Vert \mathcal{F}(z) \Vert_{B(l^2)}^2 \, \Vert \underline h \Vert^2_{\underset{1}{\overset{\infty}{\oplus}}\, H^2(\Omega)}. \notag
\end{align}

\medskip
\noindent For  $\underline x \in \mathrm{Ball}_1(l^2)$  and  $z \in \Omega$
\begin{align}
\left\Vert T_{\mathcal{F}}^* \left(\underline x \, \frac{k_z}{\Vert k_z \Vert_{H^2(\Omega)}}\right) \right\Vert_{\underset{1}{\overset{\infty}{\oplus}}\, H^2(\Omega)}^2 & = \left\Vert \mathcal{F}(z)^* \underline x \,
\frac{k_z}{\Vert k_z \Vert_{H^2(\Omega)}} \right\Vert^2_{\underset{1}{\overset{\infty}{\oplus}}\, H^2(\Omega)} \notag \\
& = \Vert \mathcal{F}(z)^* \underline x \Vert_{l^2}^2. \notag
\end{align}
Thus,
\begin{align}
\Vert T_{\mathcal{F}}\Vert = \Vert T_{\mathcal{F}}^* \Vert & \ge \sup_{z \in \Omega} \, \sup_{\underline x \in \mathrm{Ball}_1(l^2)} \Vert \mathcal{F}(z)^* \underline x \Vert_{l^2}^2 \notag \\
& \ge \sup_{z \in \Omega} \, \Vert \mathcal{F}(z)^*\Vert_{B(l^2)} \notag \\
& = \sup_{z \in \Omega} \, \Vert \mathcal{F}(z)\Vert. \notag
\end{align}
\vskip-1em
\end{proof}

\medskip
For a Hilbert space, $K$, and vectors $x,y,h \in K$, we let $x \otimes y$ denote the rank one operator defined on $K$ by
\[
(x \otimes y)(h)= \langle h \,,y\rangle x.
\]
The next lemma will be used repeatedly with $A= T_F^H $ and $k=H, \; \mbox{for} \;H \in \mathcal{H}$.
\begin{lemma}
Assume that for  $A \in B(K)$ and $k \in K \; \mbox{with} \; \Vert k \Vert_K=1$, $A \, A^* \ge \delta^2 k \otimes k$.  Then  there exists  $u_k \in (\mathrm{Ker} \, A)^{\perp}$,  so that  $A \,  u_k = k$  and  $\Vert
u_k \Vert_K \le \frac 1{\delta} $.
\end{lemma}

\begin{proof}
By the Douglas Range Inclusion Theorem, see \cite{D}, there exists a $ C \in B(H,\mathrm{Ker} \, A^{\perp})$  such that  $A \, C = k \otimes k$  and  $ \Vert C \Vert \le \frac 1{\delta}$.  Let  $ u_k = C k$.
\end{proof}

\medskip
\begin{lemma}
For $f$ a positive, bounded, lower semi-continuous function on $\partial \Omega$, there exists a nonvanishing $H \in H^{\infty}(\Omega) $, so that
\[
f = |H|^2 \; \, \sigma \mbox{-a.e.} \; \mbox{on} \; \partial \Omega.
\]
\end{lemma}
\begin{proof}
For  $\Omega = D^n$  and  $ \partial \Omega = T^n$, this is a result of Rudin \cite{R1}.  For  $ \Omega = B^n$  and  $\partial \Omega = \partial B^n$, this is a theorem of Alexandrov (see Rudin \cite{R2}, p. 32).
\end{proof}

\medskip
Recall that  \; $\mathcal{H} \triangleq \{ H \in H^{\infty}(\Omega) : \, H$ nonvanishing in $\Omega$, $ \frac 1{H} \in L^{\infty}(\partial \Omega, d \sigma)$, \, and \, $\Vert H \Vert_2 = 1 \}$.

\medskip
\noindent
For  $\{ a_j\}_{j=1}^{\infty}$  a fixed countable dense set in $\Omega$ with  $a_1 = 0$, define for each  $N = 1,2, \dots$
\[
\mathcal{C}_N \triangleq c \, o \, \{ \frac{|k_{a_j}|^2}{\Vert k_{a_j}\Vert_2^2} \, : \, j = 1, \dots, N\}
\]
Here  $k_a(\cdot)$  is the reproducing kernel for  $H^2( \Omega)$.  It is clear that $\mathcal{C}_N$ is compact and convex in $L^1(\partial \Omega,d\sigma)$.

\medskip
\noindent
Calculating, we see that for  $\Omega = D^n$   and  $g \in \mathcal{C}_N$, we have
\[
0 < \left(\frac{1 - \Vert a \Vert}{1 + \Vert a \Vert}\right)^n \le g(z) \le \left(\frac{1 + \Vert a \Vert}{1 - \Vert a \Vert}\right)^n < \infty \; \, \mbox{for all } z \in \Omega,
\]
where \; $ \Vert a \Vert = \max \, \{ \Vert a_j\Vert \, : \, j = 1, \dots, N \}$.

\medskip
For   $\Omega = B^n$  and  $g \in \mathcal{C}_N$, we have
\[
0 < \left(\frac{1 - \Vert a \Vert}{1 + \Vert a \Vert}\right)^n \le g(z) \le \left(\frac{1 + \Vert a \Vert}{1 - \Vert a \Vert}\right)^n < \infty \; \, \mbox{for all } z \in \Omega,
\]
where \; $ \Vert a \Vert = (\sum_{j=1}^{N} \,  \Vert a_j\Vert_2^2)^\frac{1}{2} $.

\medskip
\noindent
Note that for $g \in \mathcal{C}_N$, the above calculation shows that, as sets, $P^2(g\, d\sigma) \; \mbox{equals} \; H^2(\Omega)$.

Assume that  $T_F T_F^* \ge \delta^2 1 \otimes 1$  and choose  $\underline x \in \underset{1}{\overset{\infty}{\oplus}} \, H^2(\Omega)$

\noindent so that $T_F \, \underline x = 1$ and $\Vert \underline x\Vert_2 \le \frac 1{\delta}$.

\medskip
\noindent
For $  N = 1,2, \dots$ \; define
\[
\mathcal{F}_N \, : \, \mathcal{C}_N \times \underset{1}{\overset{\infty}{\oplus}} \, H^2(\Omega) \rightarrow [0, \infty)
\]
by
\[
\mathcal{F}_N(g, \underline a) \triangleq \int_{\partial \Omega} \, \Vert \underline x - P_{\mathrm{Ker}(T_F)} \underline a\Vert_{l^2}^2 \, g \, d\sigma
\]
for  $ g \in \mathcal{C}_N$  and  $\underline a \in \underset{1}{\overset{\infty}{\oplus}} \, H^2(\Omega)$.

\medskip
\noindent
Since  $g \in \mathcal{C}_N$,
\[
\underline x - P_{\mathrm{Ker}(T_F)} \underline a \in \underset{1}{\overset{\infty}{\oplus}} \, H^2(\Omega) \,\; \mbox{and } \,\, \mathcal{F}_N(g, \underline a) \, \, \mbox{is finite and positive}.
\]
For fixed  $ \underline a \in \underset{1}{\overset{\infty}{\oplus}} \, L^2(d\sigma), \, g \mapsto \mathcal{F}_N(g, \underline a)$  is linear and thus concave on the compact convex set $\mathcal{C}_N$.  For fixed  $ g \in \mathcal{C}_N$, $\underline a \mapsto \mathcal{F}(g, \underline a)$  is convex and continuous on  $\underset{1}{\overset{\infty}{\oplus}} \, H^2(\Omega)$.

\medskip
\begin{lemma}
Assume that  $T_F T_F^* \ge \delta^2 1 \otimes 1$.  For each  $N = 1, 2, \dots$,
\[
\underset{{\underline a \in \underset{1}{\overset{\infty}{\oplus}} \, H^2(\Omega)}}{\mathrm{inf}} \, \sup_{g \in \mathcal{C}_N} \, \mathcal{F}_N(g, \underline a) = \sup_{g \in \mathcal{C}_N} \, \underset{{\underline a \in \underset{1}{\overset{\infty}{\oplus}} \, H^2(\Omega)}}{\mathrm{inf}} \, \mathcal{F}_N(g, \underline a).
\]
\end{lemma}
\begin{proof}
By our remarks above, we may apply von Neumann's minimax theorem.  See, for example, Gamelin \cite{G}.
\end{proof}

\medskip
We are now ready to present our extension of Amar's theorem.
\begin{theorem}
Assume that for some $\delta >0$,  $T_{F}^H ( T_{F}^{H})^* \ge \delta^2 H \otimes H$ \; for all   $H \in \mathcal{H}$.  Then there exists a  $ G \in \underset{1}{\overset{\infty}{\oplus}} \, H^2(\Omega)$ \, with
\[
F \, G \equiv 1 \; \mbox{in } \, \Omega \; \, \mbox{and} \;\, \sup_{z \in \Omega} \Vert G(z)\Vert_{l^2} \le \frac 1{\delta}.
\]
That is,
\[
T_F  T_{G^T} \equiv I \; \; \mbox{in} \; H^2(\Omega).
\]
\end{theorem}

\begin{proof}
Since $T_F  T_F^* \ge \delta^2 1 \otimes 1$, we may choose $\underline x \in \underset{1}{\overset{\infty}{\oplus}} \, H^2(\Omega)$ so that $T_F \, \underline x = 1$ and $\Vert \underline x_0\Vert_2 \le \frac 1{\delta}$.  \\

Fix any positive integer, $N$, and any $g \in \mathcal{C}_N$. By Lemma 3, we may find an $H \in \mathcal{H}$, so that $|H|^2 = g$ \, $\sigma$-a.e. on $\partial \Omega$.

\medskip
\noindent
By our assumption
\[
T_{F}^H (T_{F}^{H})^* \ge \delta^2 H \otimes H,
\]
so there exists an $\underline x_{H} \in \underset{1}{\overset{\infty}{\oplus}} \, H^2(\Omega)$ with
\begin{equation}
T_{F}^H  (\underline x_{H}) = 1 \;\; \mbox{and} \;\; \Vert \underline x_{H}\Vert_{2, g \, d\sigma} \le \frac 1{\delta}.
\end{equation}

\medskip
Since $\underline x- \underline x_{H} \, \in \mbox{Ker}(T_F)$, we have $\underline x_{H} - \underline x = P_{\mathrm{Ker}(T_F)} \, \underline \alpha \,$ for $\, \underline \alpha=\underline x- \underline x_{H}$.
Thus (1) says that
\[
\int_{\partial \Omega} \, \Vert \underline x - P_{\mathrm{Ker}(T_F)} \, \underline \alpha \Vert^2 g \, d\sigma = \mathcal{F}_N(g, \underline \alpha) \le \frac 1{\delta^2}.
\]
Since this is true for every $g \in \mathcal{C}_N$, we may apply the minimax theorem, Lemma 4, and deduce that
\begin{equation}
\underset{{\underline \alpha \in \underset{1}{\overset{\infty}{\oplus}} \, H^2(\Omega)}}{\mathrm{inf}} \, \sup_{g \in C_N} \, \mathcal{F}_N(g, \underline \alpha) \le \frac 1{\delta^2}.
\end{equation}
Then using (2), choose  $ \underline \alpha \in \underset{1}{\overset{\infty}{\oplus}} \, H^2(\Omega)$  so that
\begin{equation}
\int_{\partial \Omega} \, \Vert \underline x - P_{\mathrm{Ker}(T_F)} \, \underline \alpha \Vert_{l^2}^2 \, g \, d \sigma \le \left( \frac 1{\delta^2} + \frac 1{N}
\right) \; \; \mbox{for all } \, g \in \mathcal{C}_N.
\end{equation}
Since  $\frac{|k_{a_j}|^2}{\Vert k_{a_j}\Vert_2^2} \in \mathcal{C}_N$  for $j = 1,2, \dots, N$, we see that

\noindent if
\[
G^{(N)} \triangleq \underline x - P_{\mathrm{Ker}(T_F)} \, \underline \alpha \, ,
\]
then
\begin{align}
 \mbox{(a)}& \; \Vert G^{(N)} (a_j)\Vert_{l^2}^2  \le \int_{\partial \Omega} \Vert G^{(N)}\Vert_{l^2}^2 \frac{|k_{a_j}|^2}{\Vert k_{a_j}\Vert_2^2} \, d\sigma \le \frac 1{\delta^2} + \frac 1{N},\, \mbox{for} \, j = 1,2, \dots, N \notag \\
 \mbox{(b)}&  \qquad \qquad \Vert G^{(N)}\Vert_2^2  \le \frac 1{\delta^2} + \frac 1{N},\; \; \mbox{and} \notag \\
 \mbox{(c)}&  \qquad \qquad F \, G^{(N)}  \equiv 1 \; \mbox{in } \, \Omega. \notag
\end{align}
Repeating this argument for each $N = 1,2,\dots$, we get a sequence of elements, $G^{(N)} \in \underset{1}{\overset{\infty}{\oplus}} \, H^2(\Omega)$, satisfying (a), (b), and (c).

\medskip
By relabeling the sequence of elements, $ \{G^{(N)} \}$, if necessary, let  $G$  be a weak  limit  of  $\{ G^{(N)} \}_{N=1}^{\infty}$  in $\underset{1}{\overset{\infty}{\oplus}} \, H^2(\Omega)$.  Fix any  $a_p \in \{a_j\}_{j=1}^{\infty}$.  Then
\[
\Vert G (a_p)\Vert_{l^2} = \lim_{N \rightarrow \infty} \, \Vert G^{(N)} (a_j) \Vert_{l^2} \le \frac 1{\delta} \;\; \mbox{by (b)}.
\]
Since  $G$  is continuous in $\Omega$  and  $\{a_j\}_{j=1}^{\infty}$  is dense in $\Omega$, we have shown that
\[
\sup_{z \in \Omega} \, \Vert G (z) \Vert_{l^2} \le \frac 1{\delta}.
\]

 \medskip
 \noindent
By (c),
\[
I= \underset{N\rightarrow\infty}{\wlim} \,T_F(G^{(N)})=T_F(G).
\]
Thus, by Lemma 1, $T_F T_G=I$.
  This completes the proof of Theorem 1.
\end{proof}

For the next theorem, we need the fact that  $ \mathrm{Ker} \, T_F = \mathrm{Ran} \, T_{\mathcal{F}}$, for an appropriate  analytic $\mathcal{F}$ . For  $\Omega = B^n$, the unit ball in $\mathbb{C}^n$, the fact that  $\mathrm{Ker} \, T_F = \mathrm{Ran} \, T_{\mathcal{F}}$   follows from results of Andersson and Carlsson \cite{AC}. For  $\Omega = D^2$,  $ \mathrm{Ker} \, T_F = \mathrm{Ran} \, T_{\mathcal{F}} $ follows from Taylor spectrum results of Putinar \cite{Pu}. That $ \mathrm{Ker} \, T_F = \mathrm{Ran} \, T_{\mathcal{F}} $ in the general case, $\Omega = D^n$, follows from an extension of the techniques of Trent \cite{Tr} and will appear in a forthcoming paper concerning the Taylor spectrum of $T_F$.

\medskip
The following shows that the Corona theorem for the polydisk or unit ball, reduces to an estimation of a lower bound for  $T_F^H (T_F^{H})^*$  where  $ H \in \mathcal{H}$  , but $H$ is not cyclic for $H^{2}(\Omega)$. (Note that we always have  $\frac 1{H} \in L^{\infty}(\partial \Omega, d\sigma)$.)
\begin{theorem}
For  $H \in \mathcal{H}$  and  $H$ cyclic in $ H^2(\Omega)$, then
\[
T_F T_F^* \ge \delta^2 1 \otimes 1 \, \Rightarrow \, T_F^H (T_F^{H})^* \ge \delta^2 H \otimes H.
\]
\end{theorem}
\begin{proof}
To show that, when $H$ is cyclic,  $T_F^H (T_F^H)^* \, \ge \, \delta^2 H \otimes H$, it suffices to find a $\underline u_H \in \underset{1}{\overset{\infty}{\oplus}} \, H^2(\Omega)$, satisfying
\begin{align}
& F\underline u_H = 1 \quad (\mbox {so} \; F(H\underline u_H) = H)  \notag \\
\mbox{and } \quad &  \, \Vert H\underline u_H \Vert_{\underset{1}{\overset{\infty}{\oplus}} \, H^2(\Omega)} \le \frac 1{\delta} \, \Vert H   \Vert_{H^2(\Omega)}= \frac{1}{\delta}.
\end{align}

Let $\underline x \, = \, T_F(T_F T_F^*)^{-1}1$.  Then such a $\underline u_H$ must have the form  $\underline u_H = \underline x - P_{\mathrm{Ker}(T_F)} \underline{\alpha}$  for some  $\underline{\alpha} \in \underset{1}{\overset{\infty}{\oplus}} \, H^2(\Omega)$.  To see that such an  $\underline{\alpha}$ exists, satisfying (4), we compute
{\allowdisplaybreaks
\begin{align}
& \underset{{\underline{\alpha} \in \underset{1}{\overset{\infty}{\oplus}} \, H^2(\Omega)}}{\mathrm{inf}}\int_{\partial \Omega} \, \Vert \underline x - P_{\mathrm{Ker}(T_F)}  \underline{\alpha} \Vert_{l^2}^2 \, |H|^2 \, d\sigma \notag \\
& \; = \, \underset{{\underline \alpha \in \underset{1}{\overset{\infty}{\oplus}} \, H^2(\Omega)}}{\mathrm{inf}}\int_{\partial \Omega} \, \Vert \underline x H - T_{\mathcal{F}} (H \underline \alpha ) \Vert_{l^2}^2 \,  d\sigma \;\; (\mbox{since Ker}(T_F)= \mbox{Ran}(T_{\mathcal{F}})) \notag \\
& \; = \, \underset{{\underline \beta \in \underset{1}{\overset{\infty}{\oplus}} \, H^2(\Omega)}}{\mathrm{inf}}\int_{\partial \Omega} \, \Vert \underline x H - T_{\mathcal{F}} (\underline \beta ) \Vert_{l^2}^2 \,  d\sigma \quad \mbox{(since H is cyclic)}\notag \\
& \; = \, \Vert P^{\perp}_{\mbox{ran}(T_{\mathcal{F}} )}(\underline x H)\Vert_ {\underset{1}{\overset{\infty}{\oplus}}\, H^2(\Omega)}^2 \notag \\
& \; = \, \Vert P^{\perp}_{\mbox{ker}(T_F )}(\underline x H)\Vert_ {\underset{1}{\overset{\infty}{\oplus}}\, H^2(\Omega)}^2 \quad \mbox{(since} \,  \mbox{Ran}(T_{\mathcal{F}})\, = \,\mbox{Ker}(T_F)) \notag \\
& \; = \, \Vert P_{\mbox{ran}(T_F^* )}(H\underline x)\Vert_ {\underset{1}{\overset{\infty}{\oplus}}\, H^2(\Omega)}^2 \notag \\
& \; = \, \Vert T_F^*(T_F T_F^*)^{-1}T_F H T_F^*(T_F T_F^*)^{-1} \underline 1\Vert_{\underset{1}{\overset{\infty}{\oplus}}\, H^2(\Omega)}^2 \notag \\
& \; = \, \Vert T_F^*(T_F T_F^*)^{-1}H\Vert_{\underset{1}{\overset{\infty}{\oplus}}\, H^2(\Omega)}^2 \notag \\
& \; \le \frac 1{{\delta}^2} \, \Vert H   \Vert^2_{H^2(\Omega)}= \frac{1}{\delta^2}.\notag
\end{align}
}
\vskip-1em
\end{proof}

\medskip
In the case that $n=1$, we may choose  $H$ in Lemma 3 to be outer and thus cyclic for $H^2(D)$. So Carleson's corona theorem for $H^{\infty}(D)$ follows from Theorems 1 and 2.

\bigskip
A very natural and interesting question arises from our work.  Thanks to Treil's remarkable example \cite{T}, we know that for an analytic \\
\, $ \mathcal{F} = [f_{ij}]_{i,j = 1}^{\infty}$  with
\[
\epsilon^2 I_{l^2} \le \mathcal{F}(z) \mathcal{F}(z)^* \le I_{l^2} \;\; \mbox{for all } \, z \in \Omega,
\]
there does not necessarily exist an analytic  $\mathcal{G} = [ g_{ij}]_{i,j = 1}^{\infty}$
\begin{align}
\mbox{with}  \quad & \mathcal{F}(z) \mathcal{G}(z) = I_{l^2} \quad \; \mbox{for all } \, z \in \Omega \notag \\
\mbox{and}  \quad & \sup_{z \in \Omega} \, \Vert \mathcal{G}(z) \Vert _{B(l^2)} < \infty. \notag
\end{align}

\medskip
How do we know when such a  $\mathcal{G}$  must exist?  For the case of the unit disk, $D$,  it is necessary and sufficient that there exist a  $\delta > 0$  with
\[
\delta^2 I \le T_{\mathcal{F}} T_{\mathcal{F}}^*.
\]

\medskip
For the polydisk and ball in $ \mathbb{C}^n$, a natural question is:  Does  $T_{\mathcal{F}}^H T_{\mathcal{F}}^{H*} \ge \delta^2I_H$ for some $\delta >0$ and for all  $H \in \mathcal{H}$  imply the existence of a bounded analytic Toeplitz operator  $T_{\mathcal{G}}$  with
\[
T_{\mathcal{F}} T_{\mathcal{G}} = I_{l^2}?
\]

\medskip
For  $\mathcal{T}(z)$, a $q \times \infty$ matrix with $q < \infty$, a modification of our techniques works, but we only get an estimate
\[
\Vert T_{\mathcal{G}}\Vert \le \frac q{\delta}.
\]

\end{document}